# REGULARITY OF LIMITS OF NONCOLLAPSING SEQUENCES OF MANIFOLDS

VITALI KAPOVITCH

ABSTRACT. We prove that iterated spaces of directions of a limit of a non-collapsing sequence of manifolds with lower curvature bound are topologically spheres. As an application we show that for any finite dimensional Alexandrov space $X^n$ with $n \geq 5$ there exists an Alexandrov space $Y$ homeomorphic to $X$ which can not be obtained as such a limit.

## 1. INTRODUCTION

The study of Alexandrov spaces with curvature bounded from below while being interesting as a subject in itself, has produced a lot of applications to classical Riemannian geometry. ( see [BGP92] for the basics of the theory of Alexandrov spaces). One of the major sources of these applications is provided by the combination of the following by now well-known facts [BGP92]:

1. Let $\mathcal{M}_k^{n,D}$ be the class of n-dimensional Riemannian manifolds with sectional curvatures bounded below by $k$ and diameter $\leq D$. Then this class is precompact in the Gromov-Hausdorff topology.
2. The property of a metric space to have curvature bounded from below is stable under taking Gromov-Hausdorff limits.

Suppose we have a sequence of manifolds $M_i^n \in \mathcal{M}_k^{n,D}$ converging to a boundary space $X$. It is not hard to see that the Hausdorff dimension of $X$ can not be greater than $n$. If it is equal to $n$ we say that the sequence $M_i^n$ converges without collapse and if it is less than $n$ we say that this sequence collapses.

The first case is understood fairly well at least topologically due to the stability theorem by Perelman [Per91], which says that for sufficiently large indices, Hausdorff approximations $M_i \to X$ are close to homeomorphisms. (Moreover, Perelman proved that these homeomorphisms can be chosen to be bi-Lipschitz but the proof of this statement has never been published.)

This result immediately implies the following finiteness theorem due to Grove-Petersen-Wu which was originally proved by other means.

**Theorem 1.1.** [GPW91] *The class of n-dimensional Riemannian manifolds with sectional curvature $\geq k$, diameter $\leq D$ and volume $\geq v$ has only finitely many topological (differentiable if $n \neq 4$) types of manifolds.*

Let $\mathfrak{M}_k^n$ denote the class of all compact n-dimensional Alexandrov spaces which can be obtained as limits of compact n-dimensional Riemannian manifolds with sectional curvatures bounded below by $k$. Let $\mathfrak{S}_k^n$ be the subclass of $\mathfrak{M}_k^n$ consisting of spaces that can be approximated by standard n-spheres with Riemannian metrics with curvatures bounded below by $k$.

**Definition 1.2.** *Let $X^n$ be an Alexandrov space of curv $\geq k$. We say that the metric on $X^n$ is smoothable if $X$ belongs to $\mathfrak{M}_{k'}^n$ for some real number $k'$.*

Notice that according to Perelman's stability theorem a smoothable Alexandrov space is a topological manifold.

Peterson conjectured in [Pet] that for any point $p$ in a smoothable space $X$ the space of directions at $p$ is homeomorphic to a sphere. Note that this does not follow from the Perelman's result since there exist Alexandrov spaces that are topological manifolds with spaces of directions at some points different from the spheres ( see example 1.5 below).

In this paper we prove the following strengthened version of Petersen's conjecture:

**Theorem 1.3.** *Suppose $X^n$ is a smoothable Alexandrov space. Then for any $x_0 \in X^n$ the space of directions $\Sigma_{x_0}X$ belongs to $\mathfrak{S}_{-1}^{n-1}$. In particular $\Sigma_{x_0}X$ is homeomorphic to a sphere of dimension $n-1$.*

Notice that since $\mathfrak{S}_{-1}^{n-1} \subset \mathfrak{M}_{-1}^{n-1}$, an obvious induction immediately yields

**Corollary 1.4.** *Suppose $X^n$ is a smoothable Alexandrov space. Then for any $x_0 \in X^n, x_1 \in \Sigma_{x_0}X, \ldots, x_i \in \Sigma_{x_{i-1}} \ldots \Sigma_{x_0}X$ the iterated space of directions $\Sigma_{x_i}\Sigma_{x_{i-1}} \ldots \Sigma_{x_0}X$ belongs to $\mathfrak{S}_{-1}^{n-i-1}$. In particular every iterated space of direction is a topological sphere.*

It is not hard to construct examples of Alexandrov metric on $\mathbb{S}^n$ that do not satisfy the conclusion of 1.4 provided $n \geq 5$.

**Example 1.5.** Let $\Sigma^3$ be the Poincare homology sphere with the metric of constant curvature 1. By taking multiple spherical suspensions of this metric space we obtain metrics of curvature $\geq 1$ on $S^m\Sigma^3$ for any $m \geq 1$. On the other hand $\Sigma^3$ is a homology sphere and hence, by the result of Edwards (cf. [Dav86]) for $m \geq 2$ the space $S^m\Sigma^3$ is known to be homeomorphic to $S^{m+3}$. Corollary 1.4 implies that the above metric on $S^m\Sigma^3$ is nonsmoothable.

We use metrics constructed in 1.5 to show that in fact nonsmoothable metrics are fairly common:

**Corollary 1.6.** *For every Alexandrov space $X^n$ with $n \geq 5$ there exists a nonsmoothable Alexandrov space $Y^n$ which is homeomorphic to $X$.*

**Remark 1.7.** In light of this Corollary it would be interesting to know how generic nonsmoothable metrics are among all Alexandrov metrics on $X$. To this end we show (see Corollary 5.3 below) that if $n = 4m + 1 \geq 5$ then for any Alexandrov space $X^n$ there exists a nonsmoothable Alexandrov space $Y$ which is homeomorphic to $X$ and such that for any fixed $k \in R$ all $n$-dimensional spaces of $curv \geq k$ sufficiently close to $Y$ are also nonsmoothable. On the other hand, the author suspects that if $X$ is homeomorphic to a smooth manifold nonsmoothable metrics are *not* dense among all Alexandrov metrics on $X$ with a fixed lower curvature bound. For example, it seems likely that every Alexandrov space of curvature $\geq 1$ sufficiently Gromov-Hausdorff close to the round sphere of constant curvature 1 is smoothable. Here is an indication why this might be true. An easy volume comparison argument shows that all iterated spaces of directions for such space have volumes almost equal to the volume of the round spheres of curvature 1 of appropriate dimension and therefore are topologically spheres. Thus no examples of the kind used in the proof of Corollary 1.6 are possible here.

Theorem 1.3 naturally leads to the following

**Question 1.8.** *Does the converse to Theorem 1.3 hold? In other words, suppose $M^n$ is an Alexandrov space such that for any point $x \in M$ the space of directions $\Sigma_{x_0} X$ belongs to $\mathfrak{S}^{n-1}_{-1}$. Is it true that the metric on $M$ is smoothable?*

One can ask an even more ambitious question:

**Question 1.9.** *Does any $n$-dimensional Alexandrov manifold with all iterated spaces of directions homeomorphic to spheres belong to $\mathfrak{M}^n_k$?*

The author suspects that this is probably false but certainly no counterexamples are known at this point.

**Remark 1.10.** Note that Theorem 1.3 might possibly give a little more than just topological restrictions on the spaces of directions of smoothable Alexandrov spaces. Namely, suppose there exists an exotic sphere $\Sigma^{n-1}$ with a Riemannian metric of $sec \geq 1$. Then its spherical suspension $S\Sigma^{n-1}$ is not smoothable and thus provides a counterexample to 1.9.

Indeed, suppose $S\Sigma^{n-1}$ is smoothable. Then the space of directions at any of the two cone points is isometric to $\Sigma^{n-1}$, hence, by Theorem 1.3 we have $\Sigma^{n-1} \in \mathfrak{S}^{n-1}_{-1}$. But this means that $\Sigma^{n-1}$ is a limit of standard $(n-1)$-spheres with smooth metrics of sectional curvature bounded below and hence, by the result of Yamagucci [Yam91], this implies that $\Sigma^{n-1}$ is diffeomorphic to a standard sphere.

In view of Corollary 1.6 one would also like to know what happens in dimensions 2, 3 and 4. By an old result of Alexandrov [AZ67], every Alexandrov space of dimension 2 is smoothable. In dimension $n = 3$ it is reasonable to expect that every Alexandrov manifold is smoothable as well. However, for $n = 4$ the situation is much less clear since it is not even known whether every four-dimensional manifold admits an Alexandrov metric of curvature bounded from below. This is most likely false because as was observed by Perelman, a positive answer would provide a counterexample to the three-dimensional Poincare conjecture.

Indeed, Perelman's stability theorem implies that a small neighborhood of a point in an Alexandrov space is bi-Lipschitz homeomorphic to an open ball in the tangent cone at this point. Since a Lipschitz structure on an $n$-dimensional manifold is unique for any $n \neq 4$ ([Sul79]), this easily implies that an Alexandrov space homeomorphic to $S^3$ is in fact bi-Lipschitz homeomorphic to $S^3$ with the canonical metric of constant curvature 1.

Now let $X^4$ be any 4-dimensional manifold that does not admit a Lipschitz structure (such manifolds exist according to [DS89]). Suppose $X$ admits an Alexandrov metric. It is easy to see that for any $p \in X$ its space of directions $\Sigma_p X$ is a simply connected manifold and hence is a homotopy sphere. Therefore there must exist a $p \in X$ such that $\Sigma_p X$ is not homeomorphic to $S^3$ since otherwise by above argument $X$ would admit a Lipschitz structure.

Another interesting problem is the finiteness of the number of differentiable types in Theorem 1.1 for $n = 4$. Differentiable finiteness for $n \neq 4$ in Theorem 1.1 is derived from topological finiteness using the fact that a compact topological manifold of dimension $\neq 4$ admits only finitely many smooth structures (cf. [KS77]). Theorem 1.3 seems to indicate that a more direct geometric argument might help to prove differentiable finiteness or even differentiable stability for $n = 4$.

In this regard let us also mention that it is currently not known if there exist Lipschitz 4-manifolds that admit nonequivalent smooth structures. If no such manifolds exist then Perelman's stability theorem implies differentiable finiteness for $n = 4$.

The author is grateful to Karsten Grove and Luis Guijarro for many helpful conversations during the preparation of this paper.

## 2. Approximation results

Let us establish some basic properties of the class $\mathfrak{S}_{-1}^{n-1}$. First of all, observe that from its definition it is obvious that $\mathfrak{S}_{-1}^{n-1}$ is closed with respect to the Gromov-Hausdorff topology among all $n$-dimensional Alexandrov spaces of $curv \geq -1$, i.e if $X^n$ is a Gromov-Hausdorff limit of spaces from $\mathfrak{S}_{-1}^{n-1}$ then $X \in \mathfrak{S}_{-1}^{n-1}$. We will use this simple observation repeatedly throughout the rest of this paper.

A basic example of a space from $\mathfrak{S}_{-1}^{n-1}$ is given by the following Lemma (cf. [Buj76]):

**Lemma 2.1.** *Let $M^n$ be a complete Riemannian manifold of $sec \geq -1$. Let $f \colon U \to \mathbb{R}$ be a proper strictly convex function on a domain $U \subset M$. Then for any $t \neq \min\limits_{U} f$ the level set $\{f = t\}$ with the induced inner metric belongs to $\mathfrak{S}_{-1}^{n-1}$.*

The proof of 2.1 is simple modulo the following two technical results that we will state here in detail for reader's convenience.

**Theorem 2.2.** [GW75] *Let $M^n$ be a Riemannian manifold and suppose $f \colon U \to \mathbb{R}$ is a proper $\delta$-strictly convex function on some domain $U \subset M$. Then there exists a sequence of smooth functions $f_m \colon U \to \mathbb{R}$ such that for any compact subset $K \subset U$ we have*

*1. $f_m$ is strictly $\delta/2$- convex on $K$ for large $m$, and*
*2. $f_m \underset{m \to \infty}{\longrightarrow} f$ uniformly on $K$.*

Another technical tool is the following special case of the metric convergence theorem by Petrunin ([Pet97, Theorem 1.2])

**Theorem 2.3.** *Let $X_m^n \xrightarrow[m \to \infty]{G-H} X^n$ be a convergent sequence of Alexandrov spaces with boundaries. Then $\partial X_m^n \xrightarrow[m \to \infty]{G-H} \partial X^n$ with respect to the induced inner metrics.*

**Remark 2.4.** Note that it is still unknown whether a level set of a convex function on an Alexandrov space taken with the induced inner metric is again an Alexandrov space. Thus the level sets in the above theorem in general are not known to be Alexandrov spaces.

*Proof of 2.1.* By Theorem 2.2 we can find a sequence of strictly convex smooth functions $f_m$ uniformly converging to $f$ on compact sets. Fix a $t$ satisfying conditions of 2.1. Clearly $\{f_m \leq t\}$ is a sequence of $n$-dimensional Alexandrov spaces of $curv \geq -1$ converging to the $n$-dimensional Alexandrov space $\{f \leq t\}$. By Petrunin's Theorem 2.3 their boundaries taken with the induced inner metrics also converge i.e

$$\{f_m = t\} \xrightarrow[m \to \infty]{G-H} \{f = t\}$$

By Gauss's formula, the level sets $\{f_m = t\}$ are smooth manifolds of sectional curvature $\geq -1$. As level sets of proper strictly convex smooth functions they are obviously diffeomorphic to $\mathbb{S}^{n-1}$ and therefore $\{f = t\} \in \mathfrak{S}^{n-1}_{-1}$.     □

## 3. CONCAVITY OF DISTANCE FUNCTIONS ON ALEXANDROV SPACES

In [Per93] Perelman introduced the following definition

**Definition 3.1.** *A function $f \colon U \to \mathbb{R}$ defined on a domain $U$ in an Alexandrov space $X$ is called $\lambda$-concave if for any unit speed shortest geodesic $\gamma \subset U$ the function $t \mapsto f(\gamma(t)) + \lambda t^2$ is concave.*

**Example 3.2.** A basic example of $\lambda$-concave functions on an Alexandrov space is given by the distance functions. Indeed. Toponogov triangle comparison implies that distance functions in a space of curvature $\geq k$ are more concave than distance functions in the model space of constant curvature $k$ and therefore it is easy to see that the following property holds:

Let $p, q \in X$ be two points in an Alexandrov space $X$ of $curv \geq k$. Let $d = d(p, q)$ and $\epsilon < d/2$. Then $f(\cdot) = d(\cdot, q)$ is $\lambda$-concave in $B(p, \epsilon)$ where $\lambda$ depends *only* on $d$ and the lower curvature bound $k$.

Here are some obvious basic properties of $\lambda$-concave functions.

- Just like for concave functions, a positive linear combination and an infimum of a family of $\lambda$-concave functions are again $\lambda$-concave
- A pointwise limit of $\lambda$-concave functions is $\lambda$-concave

**Example 3.3.** The class of examples of $\lambda$-concave functions given by Example 3.2 can be enlarged using the following simple observation from [Per93]: If $f$ is $\lambda$-concave and $\phi \colon \mathbb{R} \to \mathbb{R}_+$ is a concave $C^2$ function satisfying $0 \leq \phi' \leq 1$ then $\phi(f)$ is again $\lambda$-concave. Indeed, it is clearly enough to consider $f \colon \mathbb{R} \to \mathbb{R}$. If $f$ is $C^2$ then $\lambda$-concavity of $f$ is equivalent to the inequality $f'' \leq -\lambda$. Computing the second derivative of $\phi(f)$ we observe:

$$\phi(f)'' = \phi''(f)(f')^2 + \phi'(f)f'' \leq f'' \leq -\lambda$$

The general case immediately follows from this one since any $\lambda$-concave function on $\mathbb{R}$ can be approximated by $C^\infty$ $\lambda$-concave functions.

## 4. Proof of Theorem 1.3

The goal of this chapter is to prove Theorem 1.3 stated in the Introduction. Let us describe the strategy of the proof.

First of all, we can quickly reduce the situation to the case of pointed convergence of Riemannian manifolds with curvature bounded below to $T_{x_0}X$ where $x_0$ is any fixed point in the limit space $X$.

In [Per93] Perelman carried out a construction of strictly convex functions in a neighborhood of a point $p$ of a given Alexandrov space $X^n$ of $curv \geq k$ by using a special kind of averaging procedure for distance functions. This construction has a remarkable property: it is stable under Gromov-Hausdorff approximation of $X^n$ by spaces of the same dimension (cf.[PP93, Lemma4.3]). More precisely, if $Y^n$ is an Alexandrov space of $curv \geq k$ Gromov-Hausdorff close to $X^n$, then we can lift the distance functions used in the construction of $f$ and construct a function $\tilde{f}$ on Y which will be uniformly close to $f$ and strictly convex in a neighborhood of $\tilde{p}$. We use a modification of Perelman's construction to obtain a sequence of strictly convex functions $f_m$ on $T_{x_0}X$ satisfying the following conditions:

   (i) Each $f_m$ can be lifted to a strictly convex function $\tilde{f}_m^i$ on $M_i^n$ .
   (ii) Properly rescaled level sets of $f_m$ converge to the space of directions $\Sigma_{x_0}X$

Now level sets of $\tilde{f}_m^i$ belong to $\mathfrak{S}_{-1}^{n-1}$ by 2.1 therefore level sets of $f_m$ belong to $\mathfrak{S}_{-1}^{n-1}$ by 2.3, and finally $\Sigma_{x_0}X$ belongs to $\mathfrak{S}_{-1}^{n-1}$ as a limit of rescaled level sets of $f_m$.

*Proof of Theorem 1.3* . Recall that we start with a sequence of manifolds $M_m^n$ with sectional curvatures bounded below by $k$ converging to an Alexandrov space $X^n$. Let $x_0 \in X^n$ be any point. Let us first of all show that we can assume that $X$ is isometric to the tangent cone $T_{x_0}X$ and the lower curvature bounds for $M_m^n$'s converge to 0.

Indeed. Denote $\epsilon_m = d_{G-H}(M_m, X)$. By assumption $\epsilon_m \underset{m\to\infty}{\longrightarrow} 0$. For each $m$ there exists a $2\epsilon_m$ -Hausdorff approximation $h_m : X \to M_m$ and a $2\epsilon_m$-Hausdorff approximation $g_m : M_m \to X$ such $h_m \circ g_m$ and $g_m \circ h_m$ are uniformly $4\epsilon_m$ close to the identity maps of $X$ and $M_m$ respectively. Let $x_m = g_m(x_0)$.

**Lemma 4.1.** *Under the above assumptions we have the following convergence*

$$(1) \qquad (\frac{1}{\sqrt{\epsilon_m}}M_m, x_m) \xrightarrow[m\to\infty]{G-H} (T_{x_0}X, o)$$

*where $o$ is the cone point of $T_{x_0}X$ .*

*Proof.* By the definition of a tangent cone

$$(\frac{1}{\sqrt{\epsilon_m}}X, x_0) \xrightarrow[m\to\infty]{G-H} (T_{x_0}X, o).$$

Denote by $d^m, d$ and $\bar{d}$ the intrinsic metrics on $M_m, X$ and $T_{x_0}X$ respectively. For any fixed $R$ we know that $\delta_m(R) = d_{G-H}(B_{\frac{d}{\sqrt{\epsilon_m}}}(x_0, R), B_{\bar{d}}(o, R)) \to 0$ as $m \to \infty$. Using the triangle inequality this yields:

$$d_{G-H}(B_{\frac{d^m}{\sqrt{\epsilon_m}}}(x_m, R), B_{\bar{d}}(o, R)) \leq \sqrt{\epsilon_m} + \delta_m(R) \to 0 \text{ as } m \to \infty$$

which proves the desired convergence (1).                                    □

¿From now on we will assume that to begin with $(M_m, x_m) \xrightarrow[m\to\infty]{G-H} (T_{x_0}X, o)$ and $sec(M_m) \geq \epsilon_m$ where $\epsilon_m \xrightarrow[m\to\infty]{} 0$.

Let us proceed with the construction of $f_m$. Fix a small $\delta > 0$ and let $\delta\prime \ll \delta$. Throughout the rest of the proof we will denote by $c(n)$ various constants depending only on $n$. We will denote by $c_i$ or $c$ various constants depending on $n, \delta$ and $X$ but not on $\delta'$.

Choose a collection $\{q_\alpha\}_{\alpha \in \mathcal{A}}$ to be a maximal $\delta$-separated net in $\Sigma_{x_0}X$. For each $q_\alpha$ choose $\{q_{\alpha\beta}\}_{\beta=1,\dots,N_\alpha}$ be a maximal $\delta'$-net in $B(q_\alpha, \delta)$. The ball here is taken in $\Sigma_{x_0}X$. Note that $N_\alpha$ can be estimated from below by

$$(2) \qquad\qquad N_\alpha \geq c(n)vol(\Sigma_{x_0}X)(\delta/\delta')^{n-1}$$

Indeed, by the Absolute Volume Comparison for Alexandrov spaces [BGP92] we have that

$$volB(q_{\alpha\beta}, \delta') \leq c(n)(\delta')^{n-1}$$

and by the Relative Volume Comparison

$$volB(q_\alpha, \delta) \geq vol(\Sigma_{x_0}X)c(n)(\delta/\pi)^{n-1}$$

But since the net is taken to be maximal, the balls $\{B(q_{\alpha\beta}, 2\delta')\}_{\beta=1,\dots,N_\alpha}$ cover $B(q_\alpha, \delta)$. Hence

$$N_\alpha \cdot c(n)(\delta')^{n-1} \geq vol(\Sigma_{x_0}X)c(n)(\delta/\pi)^{n-1}$$

and therefore

$$N_\alpha \geq c(n)vol(\Sigma_{x_0}X)(\delta/\delta')^{n-1}$$

which proves (2). We will use this estimate later in the proof.

Let $\phi_{\delta'} : \mathbb{R} \to \mathbb{R}$ be the continuous function uniquely determined by the following properties:

(1)  $\phi_{\delta'}(0) = 0$

(2) $\phi'_{\delta'}(t) = 1$ for $t \leq 1 - \delta'$
(3) $\phi'_{\delta'}(t) = 1/2$ for $t \geq 1 + \delta'$
(4) $\phi''_{\delta'}(t) = -1/(4\delta')$ for $1 - \delta' < t < 1 + \delta'$

Now define $f_{\delta'\alpha}$ by the following formula:

$$f_{\delta'\alpha}(x) = \frac{1}{N_\alpha} \sum_{\beta=1}^{N_\alpha} \phi_{\delta'}(d(x, q_{\alpha\beta}))$$

Then according to Lemma 3.6 from [Per93], (cf. [PP93, Lemma 4.3]) the function $f_{\delta'\alpha}$ is strictly $c/\delta'$-concave in $B(o, \delta'/2)$ for all sufficiently small $\delta'$ (we will reprove this statement in Lemma 4.2 below). Finally, define $f_{\delta'}$ as $f_{\delta'}(x) = \min_\alpha f_{\delta'\alpha}(x)$. Then it is clear that $f_{\delta'}$ is strictly $c/\delta'$-concave in $B(o, \delta'/2)$.

Let us examine this function more carefully. First, observe that for all $\alpha$ we get $f_{\delta'\alpha}(o) = \phi_{\delta'}(1)$ and hence $f_{\delta'}(o) = \phi_{\delta'}(1)$. Moreover, we claim that $o$ is a point of a strict local maximum of $f_{\delta'}$. Indeed, let $x \in B(o, \delta')$ be a point sufficiently close to $o$ with $d(o, x) = t$. Without too much abuse of notation we can write $x = t\xi$ for some $\xi \in \Sigma_{x_0} X$. By our construction there exists an $\alpha_0$ such that $\sphericalangle \xi q_\alpha \leq \delta$. Then clearly $\sphericalangle \xi q_{\alpha\beta} \leq 2\delta$, and therefore by the first variation formula $d(\xi, q_{\alpha\beta}) \leq 1 - t\cos(3\delta)$ for sufficiently small $t$. Hence, by monotonicity of $\phi_{\delta'}$ we immediately get

$$f_{\delta'}(t\xi) \leq \phi_{\delta'}(1 - t\cos(3\delta)) < \phi_{\delta'}(1) = f_{\delta'}(o)$$

We will give a more accurate estimate for $f_{\delta'}(t\xi)$ later.

Let us lift $f_{\delta'}$ to the elements of the sequence $(M_m, x_m)$ in a natural way. More precisely, according to (1) there exists a $\mu_m$-Hausdorff approximation $h_m : B(o, 2) \to B_{M_m}(x_m, 2)$, where $\mu_m \underset{m \to \infty}{\longrightarrow} 0$ and $h_m(o) = x_m$. Let $q_\alpha^m = h_m(q_\alpha)$ and $q_{\alpha\beta}^m = h_m(q_{\alpha\beta})$ for all $\alpha$ and $\beta$. Then we put

$$f_{\delta'\alpha}^m(y) = \frac{1}{N_\alpha} \sum_{\beta=1}^{N_\alpha} \phi_{\delta'}(d(y, q_{\alpha\beta}^m))$$

and

$$f_{\delta'}^m(y) = \min_\alpha f_{\delta'\alpha}^m(y)$$

The most important technical part in the proof of Theorem 1.3 is the following modification of Lemma 3.6 from [Per93].

**Lemma 4.2.** *For $\mu_m \ll \delta'$ we have that $f_{\delta'}^m$ is $c/\delta'$ concave in $B_{Y_m}(x_m, \delta'/2)$ where the constant $c$ is independent of $\delta'$ but it does depend on $X$ and $\delta$.*

**Remark 4.3.** Let us mention that the statement of 4.2 is essentially contained in the proof of [PP93, Lemma 4.3]. However, since the proof there is omitted and the proof of [Per93, Lemma 3.6] is not very detailed, we present a complete proof here.

Here is the informal idea of the proof. Fix a sufficiently large $m$ and let $\gamma(t)$ be a geodesic in $B_{Y_m}(x_m, \delta'/2)$ and let $u_{\alpha\beta}(t) = d(\gamma(t), q^m_{\alpha\beta})$.

Suppose for a moment that $u_{\alpha\beta}(t)$ is smooth it $t$ for all indices $\alpha\beta$. By 3.3 $(\phi_{\delta'}(u_{\alpha\beta}))(t)$ is $c_1$-concave for any $\delta' > 0$ where $c_1$ depends *only* on the lower curvature bound for $M_m$ ( recall that $d(\gamma(t), q^m_{\alpha\beta}) \approx 1$ ) and therefore $(\phi_{\delta'}(u_{\alpha\beta}))''(t) \le -c_1$.

On the other hand a volume comparison argument shows that for any fixed $t$ we have $|u'_{\alpha\beta}(t)| > c_2$ for vast majority of indices $\alpha\beta$. For all such indices we have $(\phi_{\delta'}(u_{\alpha\beta}))''(t) = \phi''_{\delta'}(u_{\alpha\beta})(u'_{\alpha\beta}(t))^2 + \phi'_{\delta'}(u_{\alpha\beta})u''_{\alpha\beta}(t) \le -\frac{1}{\delta'}(c_2)^2 + c_1$. Thus for all indices satisfying $|u'_{\alpha\beta}(t)| > c_2$ we can make $(\phi_{\delta'}(u_{\alpha\beta}))''(t)$ to be as negative as we like by taking $\delta'$ to be sufficiently small and since for the remaining indices $(\phi_{\delta'}(u_{\alpha\beta}))''(t) \le c_1$, the same is true for the average of $(\phi_{\delta'}(u_{\alpha\beta}))''(t)$.

**Remark 4.4.** Our proof will also show that $f_{\delta'}$ is strictly $c/\delta'$-concave in $B(o, \delta')$.

*Proof of Lemma 4.2:* Since a minimum of a family of concave functions is again concave it is certainly enough to prove the desired concavity property for each $f^m_{\delta'\alpha}$.

First observe the following: Let $\Gamma \subset \Sigma^{n-1}$ be any subset in a space of $curv \ge 1$. Consider the set $U_\epsilon = \{x \in \Sigma \mid \pi/2 - \epsilon \le d(x, \Gamma) \le \pi/2 + \epsilon\}$. Then $vol(U_\epsilon) \le c(n)\epsilon$. This is a direct corollary of the volume comparison for Alexandrov spaces (cf. [BGP92, Lemma 8.2]). Therefore, if $\{z_i\}_{i=1...N}$ is a maximal $\delta'$-net in $U_\epsilon$, then

$$(3) \qquad\qquad N \le \epsilon c(n) vol(\Sigma)(\delta')^{1-n}$$

Now let $xz \subset B_{Y_m}(x_m, \delta'/2)$ be a shortest curve and let $y$ be its midpoint. Let $t = d(xy)$. Consider the set of indices $I'_\alpha$ such that for any $\alpha\beta \in I'_\alpha$ we have $|\cos \sphericalangle q_{\alpha\beta}yx| > \nu$. And let $I''_\alpha$ be the set of indexes for which $|\cos \sphericalangle q_{\alpha\beta}yx| \le \nu$. Here $\sphericalangle q_{\alpha\beta}yx$ stands for the minimal possible angle between $yx$ and a shortest geodesic connecting $y$ and $q_{\alpha\beta}$.

Denote $N'_\alpha = |I'_\alpha|$ and $N''_\alpha = |I''_\alpha|$.

Note that for small $\mu_m$ and for $\beta_1 \ne \beta_2$ we certainly have $\sphericalangle q^m_{\alpha\beta_1}xq^m_{\alpha\beta_2} \ge \delta'/2$. Hence, by (3) it follows that $N''_\alpha \le \nu c(n)(\delta')^{1-n}$. On the other hand, according to (1) the total number of points $N_\alpha = N'_\alpha + N''_\alpha$ satisfies $N_\alpha \geq c(\delta')^{1-n}$. Therefore

$$(4) \qquad N''_\alpha/N_\alpha \leq \nu/c \text{ and } N'_\alpha/N_\alpha \geq 1 - \nu/c$$

Fix a small $\nu$ satifying

$$(5) \qquad \nu/c \leq 1/2$$

Now we are going to give two separate concavity estimates for $f^m_{\delta'\alpha\beta}$ along $xz$: one for $\alpha\beta \in I'_\alpha$ and the other one for $\alpha\beta \in I''_\alpha$.

First, choose any $\alpha\beta \in I'_\alpha$. In this case, we have the following estimate:

$$(6) \qquad 2f^m_{\delta'\alpha\beta}(y) - f^m_{\delta'\alpha\beta}(x) - f^m_{\delta'\alpha\beta}(z) \geq (c(\delta')^{-1}\nu^2)d(xz)^2$$

Indeed, by construction of $I'_\alpha$, we have $|\cos \sphericalangle q_{\alpha\beta} yx| > \nu$. Consider, for example, the case when $\cos \sphericalangle q_{\alpha\beta} yx > \nu$. (The other case is treated similarly by reversing the roles of $x$ and $z$). By the triangle comparison we obtain:

$$d(xq) \leq d(yq) - \nu t + c_1 t^2 \text{ and } d(zq) \leq d(yq) + \nu t + c_1 t^2$$

where $c_1$ depends only on the lower curvature bound for $M_m$ (Recall that $d(xy) \approx 1$). Next observe that by construction, $\phi_{\delta'}(t)$ is monotone for all $t$ and is strictly $\frac{1}{4\delta'}$ concave with constant second derivative for $1 - \delta' < t < 1 + \delta'$. Hence,

$$\phi_{\delta'}(d(yq) - \nu t + c_1 t^2) = \phi_{\delta'}(d(yq)) - \phi'_{\delta'}(d(yq))(\nu t - c_1 t^2) - 1/(8\delta')(\nu t - c_1 t^2)^2$$

which for sufficiently small $t$ implies

$$f^m_{\delta'\alpha\beta}(x) \leq \phi_{\delta'}(d(yq) - \nu t + c_1 t^2) \leq \phi_{\delta'}(d(yq)) - \phi'_{\delta'}(d(yq))(\nu t) - \left(c_2 \nu^2/\delta' - c_1\right) t^2$$

Similarly,

$$f^m_{\delta'\alpha\beta}(z) \leq \phi_{\delta'}(d(yq)) + \phi'_{\delta'}(d(yq))(\nu t) - \left(c_2 \nu^2/\delta' + c_1\right) t^2$$

Adding these two inequalities we immediately obtain (6). So we have a good concavity estimate for $f^m_{\delta'\alpha\beta}$ along $xz$ for any $\alpha\beta \in I'_\alpha$. For the rest of the indices 3.3 implies that $f^m_{\delta'\alpha\beta}$ is $c_3$-concave and therefore

$$(7) \qquad 2f^m_{\delta'\alpha\beta}(y) - f^m_{\delta'\alpha\beta}(x) - f^m_{\delta'\alpha\beta}(z) \geq -c_3 d(xz)^2$$

for some $c_3$ independent of $\delta'$. Combining (6) and (7) we get

$$2f^m_{\delta'\alpha}(y) - f^m_{\delta'\alpha}(x) - f^m_{\delta'\alpha}(z) \geq N'_\alpha/N_\alpha c_2(\delta')^{-1}\nu^2 d(xz)^2 - N''_\alpha/N_\alpha c_3 d(xz)^2$$

Hence, using (4) we obtain

$$2f^m_{\delta'\alpha}(y) - f^m_{\delta'\alpha}(x) - f^m_{\delta'\alpha}(z) \geq ((1 - \nu/c)c_2(\delta')^{-1}\nu^2 - \nu/cc_3)d(xz)^2$$

Finally, by (5) this last inequality implies that for $\delta' \ll \delta$ we have

$$2f^m_{\delta'\alpha}(y) - f^m_{\delta'\alpha}(x) - f^m_{\delta'\alpha}(z) \geq c_4(\delta')^{-1}d(xz)^2$$

$\square$

**Lemma 4.5.** *For any $t < \delta'/4$ the level set $(f_{\delta'} = \phi_{\delta'}(1-t))$ belongs to $\mathfrak{S}_{-1}^{n-1}$*

*Proof.* By Lemma 4.2 $f_{\delta'}^m$ is $\frac{1}{2}c\delta'$ concave in $B_{M_m}(x_m, \delta'/2)$. Recall that $sec(M_m) \geq \epsilon_m$ where $\epsilon_m \underset{m \to \infty}{\longrightarrow} 0$. Therefore $sec(M_m) \geq -1$ for large $m$ and by Lemma 2.1 level sets $f_{\delta'}^m = \phi_{\delta'}(1-t)$ belong to $\mathfrak{S}_{-1}^{n-1}$ for any $0 < t < \delta'/2$. By construction of $f_{\delta'}^m$ we obviously have that $f_{\delta'}^m \underset{m \to \infty}{\overset{\text{unif}}{\Longrightarrow}} f_{\delta'}$. By Theorem 2.3 this implies that $(f_{\delta'}^m = \phi_{\delta'}(1-t)) \underset{m \to \infty}{\overset{G-H}{\longrightarrow}} (f_{\delta'} = \phi_{\delta'}(1-t))$ where the metrics on the level sets are taken to be the induced inner metrics. Hence $(f_{\delta'} = \phi_{\delta'}(1-t)) \in \mathfrak{S}_{-1}^{n-1}$ as claimed. $\qquad \square$

**Remark 4.6.** Our proof of Lemma 4.5 actually shows that the level sets of $f_{\delta'}$ belong to $\mathfrak{S}_{-\epsilon}^{n-1}$ for any positive $\epsilon$.

**Lemma 4.7.** *The following estimate holds for all $0 < t \ll \delta'$:*

$$d_H\left(\frac{1}{t}\left\{f_{\delta'} = \phi_{\delta'}(1-t)\right\}, \Sigma_{x_0}X\right) \leq \frac{1}{\cos(3\delta)} - 1$$

*where $d_H$ means the Hausdorff distance between subsets of $T_{x_0}X$ and we identify $\Sigma_{x_0}X$ with the unit sphere at $o$ in $T_{x_0}X$ and $\frac{1}{t}\{f_{\delta'} = \phi_{\delta'}(1-t)\}$ with the homothetic image of $\{f_{\delta'} = \phi_{\delta'}(1-t)\}$ under the $1/t$ -homothety of $T_{x_0}X$.*

*Proof.* Let $x = t\xi$ for some $\xi \in \Sigma_{x_0}X$ as before. First of all observe that by the triangle inequality $d(x, q_{\alpha\beta}) \geq 1 - t$ for any $\alpha\beta$ and therefore $\phi_{\delta'}(1-t) \leq f_{\delta'}(x)$.

On the other hand, as we have seen there exists an $\alpha_0$ such that $\sphericalangle \xi q_{\alpha_0} \leq \delta$. Clearly, $\sphericalangle q_{\alpha_0\beta} \leq 2\delta$ for any $\beta$. By the cosine law we have

$$d^2(x, q_{\alpha_0\beta}) = 1 + t^2 - 2t\cos(\sphericalangle \xi q_{\alpha_0\beta}) \leq 1 + t^2 - 2t\cos(2\delta)$$

For $t \ll \delta$ this implies $d(x, q_{\alpha_0\beta}) \leq 1 - t\cos 3\delta$ and therefore

$$(8) \qquad\qquad f_{\delta'\alpha_0}(x) \leq \phi_{\delta'}(1 - t\cos 3\delta).$$

Since $f_{\delta'}(x) \leq f_{\delta'\alpha_0}(x)$ we obtain

$$(9) \qquad\qquad \phi_{\delta'}(1-t) \leq f_{\delta'}(x) \leq \phi_{\delta'}(1 - t\cos 3\delta)$$

for any $x$ with $d(x,o) = t$ and $t$-sufficiently small.

Now let $t$ be sufficiently small and let us look at the level set $\{f_{\delta'}(x) = \phi_{\delta'}(1-t)\}$. By inequality (9) we have

$$\phi_{\delta'}(1 - d(x,o)) \leq \phi_{\delta'}(1-t) \leq \phi_{\delta'}(1 - d(x,o)\cos 3\delta)$$

which by monotonicity of $\phi_{\delta'}$ implies

$$d(x,o)\cos 3\delta \leq t \leq d(x,o)$$

Therefore

$$1 \leq d(x, o)/t \leq 1/(\cos 3\delta)$$

and the conclusion of Lemma 4.7 follows. $\qquad\square$

Now we can finally finish the proof of Theorem 1.3. Let us choose a sequence of positive numbers $\delta_k \underset{k \to \infty}{\longrightarrow} 0$. Then by Lemma 4.7 we can choose $t_k \ll \delta_k'$ such that

$$d_H \left( \frac{1}{t_k} \left( f_{\delta_k'} = \phi_{\delta_k'} (1 - t_k) \right), \Sigma_{x_0} X \right) \leq \frac{1}{\cos(3\delta_k)} - 1$$

where $d_H$ stands for the Hausdorff distance between subsets of $T_{x_0} X$. Hence, by Theorem 2.3 we have that

$$\frac{1}{t_k} \left( f_{\delta_k'} = \phi_{\delta_k'} (1 - t_k) \right) \underset{k \to \infty}{\overset{G-H}{\longrightarrow}} \Sigma_{x_0} X$$

in the induced inner metrics. But by Lemma 4.5 we already know that each level set $f_{\delta_k'} = \phi_{\delta_k'} (1 - t_k)$ belongs to $\mathfrak{S}_{-1}^{n-1}$. Hence $\Sigma_{x_0} X \in \mathfrak{S}_{-1}^{n-1}$ as well. $\qquad\square$

**Remark 4.8.** Since the proof of Theorem 1.3 is local the theorem remains true in case of pointed convergence of noncompact manifolds with $\sec \geq k$ to a limit space of the same dimension.

## 5. APPLICATIONS OF THEOREM 1.3

As was mentioned in the introduction, the basic example of an Alexandrov metric on a topological manifold that does not satisfy the conditions of Corollary 1.4 is as follows. Let $\Sigma^3$ be the Poincare homology sphere. Recall that it can be constructed as a quotient of $S^3$ with the canonical metric of constant curvature 1 by a free isometric action of the icosahedral group $I^*$. Therefore it is a manifold of constant curvature 1. By taking multiple spherical suspensions of this metric space we obtain metrics of curvature $\geq 1$ on $S^m \Sigma^3$ for any $m \geq 1$. On the other hand, since $\Sigma^3$ is a homology sphere, by the result of Edwards (cf. [Dav86]) the space $S^m \Sigma^3$ is known to be homeomorphic to $S^{m+3}$ for any $m \geq 2$. Since some of the iterated spaces of directions for $S^m \Sigma^3$ are isometric to $\Sigma^3$ this space is nonsmoothable by Corollary 1.4.

This construction can be used to prove Corollary 1.6 which says that nonsmoothable metrics are fairly common.

*Proof of Corollary 1.6.* Let $X^n$ be an Alexandrov space of dimension $n \geq 5$. If $X$ is not homeomorphic to a smooth manifold the statement of Corollary 1.6 is obvious by Perelman's stability Theorem.

Now suppose that that $X$ is homeomorphic to a smooth manifold $M^n$.

Let $(\Sigma^n, d_\Sigma) = S^{n-3}\Sigma^3$ be the $n$-sphere with the nonsmoothable metric constructed above. By construction the metric on $\Sigma^n$ is smooth away from $\Sigma^{n-4} \subset \Sigma^n$. Let $x \in \Sigma^n\backslash\Sigma^{n-4}$. We can find a small neighborhood $U$ of $x$ which is a Riemannian manifold. It is clear that we can use $U$ to construct a metric on the connected sum $M\#\Sigma^n$ which is Riemannian away from $\Sigma^{n-4}$ and coincides with $d_\Sigma$ on $\Sigma^n\backslash U$. This metric is obviously Alexandrov and it still contains points with some of the iterated spaces of directions isometric to $\Sigma^3$ and thus it is nonsmoothable.

Since $M\#\Sigma^n$ is homeomorphic to $X$ the conclusion of Corollary 1.6 follows.  $\square$

The proof of Theorem 1.3 can be easily modified to prove the following

**Theorem 5.1.** *Let $M^n_m$ be a sequence of $n$-dimensional Alexandrov spaces converging without collapse to an Alexandrov space $X^n$ with $\partial X = \emptyset$. Let $x_0 \in X^n$ be any point. Then for all sufficiently large $m$ there exist points $p_m \in M^n_m$ such that $\Sigma_{x_0}X$ is homeomorphic to $\Sigma_{p_m}M^n_m$.*

To prove the statement of Theorem 5.1 we have to utilize the following general Lemma:

**Lemma 5.2.** *Let $X$ be an Alexandrov space without a boundary and let $f : U \to \mathbb{R}$ be a proper strictly concave function in some domain $U \subset X$. Let $p \in U$ be the point of strict maximum of $f$. Then $\Sigma_pX$ is homeomorphic to any nonempty level set $\{f(x) = c\}$ provided $c < f(p)$.*

*Proof.* Let us first prove Lemma 5.2 in the special case when $f$ satisfies the following extra condition:

$$(10) \qquad\qquad f(x) \leq f(p) - L \cdot d(x, p)$$

for all $x$ sufficiently close to $p$ and some fixed constant $L$.

First let us note that by Perelman's Stability Theorem we immediately obtain that different superlevel sets $\{f \geq c\}$ are homeomorphic to each other for $c < f(p)$, and since homeomorphisms of Alexandrov spaces send boundaries to boundaries (cf. [Per91, Thm 4.6]), the level sets $\{f(x) = c\}$ are homeomorphic to each other as well.

Now take any sequence of numbers $\lambda_n \to o$. Then we know that the pointed sequence of spaces $(\frac{1}{\lambda_n}X, p)$ Gromov-Hausdorff converges to $(T_pX, o)$ as $n \to \infty$. Let us denote $\frac{1}{\lambda_n}X$ by $Y_n$. Without a loss of generality we can assume that $f(p) = 0$. Consider $g_n = \frac{1}{\lambda_n}f : Y_n \to \mathbb{R}$. It is clear that this family of functions is uniformly Lipschitz and uniformly bounded in balls $B_{Y_n}(o, R)$. Therefore, by Arzela-Ascoli we can find a subsequence $g_{n_k}$ uniformly converging to $g : T_p X \to \mathbb{R}$ ( One can think of $g$ as the "differential' of $f$ ). Then, obviously, $g$ is convex and condition (10) guarantees that it has a strict maximum at $o$. Now it is clear that the radial projection along the rays emanating from $o$ will provide a homeomorphism between a nonmaximal level set of $g$ and $S(o, \epsilon)$ which is homothetic (and hence homeomorphic) to $\Sigma_p X$. On the other hand, by Perelman's Stability Theorem we know that level sets of $g$ are homeomorphic to the level sets of $g_{n_k}$ for large $k$. Now we can conclude the proof by noticing that level sets of $g_{n_k}$ are just rescaled level sets of $f$.

Next let us look at the general case when we do not assume that $f$ satisfies (10). Notice that in the proof of Lemma 4.7 we never used the fact that $X$ was a topological manifold. We can use the same construction to show that for any Alexandrov space $X$ and any point $x \in X$ there exists a Lipschitz function $h$ which is strictly concave in a neighborhood of $x$ and has its maximum at $x$ ( cf. [PP93, Lemma 4.3]). The only difference from the situation of Lemma 4.7 in constructing such a function is that we should choose points $q_\alpha$ and $q_{\alpha\beta}$ on a small metric sphere $S(x, r)$ centered at $x$. It is also easy to see that if $r$ is chosen to be sufficiently small, the constructed strictly concave function will satisfy the inequality similar to inequality (8) from the proof of Lemma 4.2 above and hence it will satisfy condition (10).

Now suppose $f$ is a strictly concave function on $U$ and let $p$ be its point of maximum. By above there exists a function $h$ which is strictly concave on a small neighborhood of $p$ such that $h$ has a maximum at $p$ and it satisfies condition (10). Consider the family of functions $f_\epsilon = f + \epsilon h$ where $\epsilon \geq 0$. It is obvious that there exists an open set $U'$ containing $p$ such that each $f_\epsilon$ is strictly concave in $U'$ and has a maximum at $p$. It is also obvious that $f_\epsilon$ satisfies condition (10) for any $\epsilon > 0$. Thus almost maximal level sets of $f_\epsilon$ are homeomorphic to $\Sigma_p X$ for any $\epsilon > 0$. On the other hand $f_\epsilon \underset{\epsilon \to 0}{\overset{unif}{\Longrightarrow}} f$ and therefore level sets of $f_\epsilon$ are homeomorphic to level sets of $f$ by Perelman's theorem. $\qquad\square$

*Proof of Theorem 5.1.* As was mentioned above the proofs of Lemma 4.1 and Lemma 4.2 never used the fact that elements of the sequence were Riemannian manifolds. Therefore, proceeding as in the proof of Theorem 1.3, we can construct families of functions $f_{\delta'}$ on $T_{x_0} X$ and their liftings $f_{\delta'}^m$ on $M_m$. Then as before we know

(1) $f_{\delta'}$ is strictly $c/\delta'$ concave in $B(o, \delta'/2)$ and has a strict maximum at $o$.

(2) For all sufficiently large $m$ we have that $f_{\delta'}^m$ is $\frac{1}{2} c/\delta'$ concave in $B_{M_m}(x_m, \delta'/2)$.

(3) $f_{\delta'}^m \underset{m \to \infty}{\overset{unif}{\Longrightarrow}} f_{\delta'}$

Perelman's Stability Theorem combined with (3) implies that superlevel sets of $f_{\delta'}$ and $f_{\delta'}^m$ are homeomorphic for all large $m$, and therefore, the level sets of $f_{\delta'}$ and $f_{\delta'}^m$ are homeomorphic for all large $m$ as well. Now let us look at superlevel sets of $f_{\delta'}$. By construction, $o$ is the point of strict maximum of $f_{\delta'}$. Hence level sets of $f_{\delta'}$ are homeomorphic to $\Sigma_{x_0}X$ by Lemma 5.2. By the same Lemma, level sets of $f_{\delta'}^m$ are homeomorphic to $\Sigma_{p_m}Y_m$ where $p_m$ is the point of strict maximum of $f_{\delta'}^m$ and by transitivity $\Sigma_{x_0}X$ is homeomorphic to $\Sigma_{p_m}M_m$. $\qquad\blacksquare$

As an immediate corollary of Theorem 5.1 we show that for $n = 4l + 1 \geq 5$ Corollary 1.6 can be strengthened in the following way:

**Corollary 5.3.** *Let $X^n$ be an Alexandrov space of dimension $n = 4l + 1 \geq 5$. Then there exists a nonsmoothable Alexandrov space $Y$ which is homeomorphic to $X$ and which satisfies the following property:*

*For any fixed $k \in \mathbb{R}$ there exists an $\epsilon > 0$ such that any Alexandrov space $Z^n$ of $curv \geq k$ with $d_{G-H}(Y, Z) \leq \epsilon$ is nonsmoothable.*

*Proof.* As before we only have to consider the case when $X$ is homeomorphic to a smooth manifold $M^n$. Consider the diagonal action of the icosahedral group $I^*$ on $\underbrace{\mathbb{R}^4 \times ... \times \mathbb{R}^4}_{l}$. This action is obviously free when restricted to $S^{4l-1}$. Let $\Sigma = S^{4l-1}/I^*$. By the same argument as in the proof of 1.6, we observe that $S^2\Sigma$ is homeomorphic to $S^n$. Again proceeding in the same way as in the proof of 1.6 we can construct an Alexandrov metric on $Y = M^n \# S^2\Sigma$ such that there is a point $p \in Y$ with $\Sigma_p Y$ isometric to $S\Sigma$.

Fix any $k \in \mathbb{R}$. By Theorem 5.1, any Alexandrov space $Z^n$ of $curv \geq k$ sufficiently close to $Y$ has a point with a space of directions homeomorphic to $S\Sigma$. Since $S\Sigma$ is obviously not homeomorphic to a sphere, $Y$ is nonsmoothable by Theorem 1.3. $\qquad\blacksquare$

**Remark 5.4.** As it was pointed out to the author by the referee, it should be possible to show that in the setting of Theorem 5.1 any iterated space of directions of $X$ is homeomorphic to an iterated space of directions of $M_m$ for large $m$. ( Notice that this would imply that Corollary 5.3 holds for any $n \geq 5$).

However, to carry out the proof one has to use a refined version of Perelman's stability theorem the proof of which goes beyond the intent of the present paper.

DEPARTMENT OF MATHEMATICS, UNIVERSITY OF PENNSYLVANIA, PHILADELPHIA, PA 19104, FAX: $(215) - 573 - 4063$

*E-mail address*: `vitali@math.upenn.edu`